\documentclass[reqno]{amsart}
\usepackage{amsmath}
\usepackage{amsfonts}

\usepackage{xy}
\usepackage[bookmarks=false,breaklinks=true]{hyperref}
\usepackage{float}
\usepackage{graphicx}
\usepackage[percent]{overpic}
\usepackage{color}

\usepackage{epigraph}
\definecolor{myurlcolor}{rgb}{0,0,0.7}
\usepackage{hyperref}
\hypersetup{colorlinks,
linkcolor=myurlcolor,
citecolor=myurlcolor,
urlcolor=myurlcolor}
\usepackage[mathscr]{euscript}

\usepackage{color}

\definecolor{lightblue}{rgb}{0.5,0.5,0.8}
\definecolor{darkblue}{rgb}{0,0,0.8}
\definecolor{darkpurple}{rgb}{0.4,0,0.5}
\definecolor{brickred}{rgb}{0.7,0,0}
\definecolor{darkred}{rgb}{0.6,0,0}
\definecolor{verydarkred}{rgb}{0.4,0,0}

\usepackage{amsfonts,amssymb}
\usepackage{amsmath}

\usepackage[english]{babel}

\newtheorem{thm}{Theorem}

\theoremstyle{definition}

\theoremstyle{remark}

\numberwithin{equation}{section}

\newcommand{\Var}{\mathsf{Var}}
\newcommand{\Mot}{\mathsf{Mot}}

\newcommand{\R}{{\mathbb R}}
\newcommand{\C}{{\mathbb C}}
\newcommand{\F}{{\mathbb F}}

\newcommand{\Z}{{\mathbb Z}}
\renewcommand{\P}{{\mathbb P}}
\renewcommand{\L}{{\mathbb L}}

\newcommand{\maps}{\colon}

\begin{document}

\title{Motivating Motives}

\author{John C.~Baez}
\address{Department of Mathematics\\ 
University of California\\ 
Riverside CA 92521\\
USA }
\email{baez@math.ucr.edu}

\begin{abstract}
Underlying the Riemann Hypothesis there is a question whose full answer still eludes us: what do the zeros of the Riemann zeta function really mean? As a step toward answering this, Andr\'e Weil proposed a series of conjectures that include a simplified version of the Riemann Hypothesis in which the meaning of the zeros becomes somewhat easier to understand. Grothendieck and others worked for decades to prove Weil's conjectures, inventing a large chunk of modern algebraic geometry in the process. This quest, still in part unfulfilled, led Grothendieck to dream of ``motives": mysterious building blocks that could explain the zeros (and poles) of Weil's analogue of the Riemann zeta function. This exposition by a complete amateur tries to sketch some of these ideas in ways that other amateurs can enjoy. 
\end{abstract}

\maketitle

\section{Introduction}

Grothendieck's ``motives'' are one of his most mysterious and tantalizing ideas.   Indeed, he felt this way himself.   In \emph{R\'ecoltes et Semaille} he wrote:
\begin{quote}
Among all the mathematical discoveries which I've been privileged to make, 
the concept of the motive
still impresses me as the most fascinating, 
the most charged with mystery---indeed at the very heart of the profound identity 
of geometry and arithmetic.
\end{quote}
One reason for this mystique is that motives arose in Grothendieck's long quest to prove the Weil Conjectures.  The most difficult of these was an analogue of the Riemann Hypothesis.   By dint of  strenuous effort, he reduced this problem to some conjectures about motives, called the ``Standard Conjectures''.  But these remain unproved even to this day!  His student Deligne finished off the Weil Conjectures by a different method.  Thus, the theory of motives remains obscure in many respects. 

Nonetheless, motives are important.    Quantum mechanics revealed
that any physical system is in a superposition of ``energy eigenstates'', each vibrating at its own frequency.  Similarly, Grothendieck discovered that you can take polynomial equations with coefficients in the field with \(p^n\) elements, where \(p\) is some prime number, and break apart their set of solutions into pieces that are not sets.  These pieces are called ``motives'': they're more like vector spaces than sets.    Each piece has something resembling a number of points---but in fact, this number is essentially the trace of some operator to the \(n\)th power.  This number can grow exponentially but also \emph{oscillate} as a function of \(n\).   So it doesn't need to be positive!   When you add up these numbers, you get the number of solutions of your equations.

My goal here is not to explain this in detail: I merely want to get beginners interested in motives.  Mathematicians have a curious turn of phrase where they speak of ``motivating'' a concept when they really mean motivating someone to study that concept.   In this language, I am trying to motivate motives.    

Since this is an explanation by an amateur for amateurs, experts should avert their eyes.  
I will ruthlessly suppress details that are not strictly necessary, while trying to not fall 
into outright error.  All references---mainly suggestions for further study---are in the last section.

\section{The Riemann Hypothesis}

The Riemann zeta function is given by a sum that converges for \(\mathrm{Re}(z) > 1\):
\[   \zeta(z) = \frac{1}{1^z} + \frac{1}{2^z} + \frac{1}{3^z} + \cdots \]
but we can analytically continue it to the whole complex plane except for a pole
at \(z = 1\).  It is zero at negative even integers, called ``trivial'' zeros.   It also has other zeros, called ``nontrivial'' zeros---and the Riemann Hypothesis says all these lie on the 
line \(\mathrm{Re}(z) = \frac{1}{2} \).  However, this bare statement sheds little light on why the Riemann Hypothesis is interesting.  

It is often claimed that there is no useful formula for the \(n\)th prime number. 
But Riemann came up with something just as good: a formula for the ``prime counting function'' \(\pi(n)\), which is the number of primes \(\le n\):

\vskip 1em

\begin{center}
\includegraphics[width=25em]{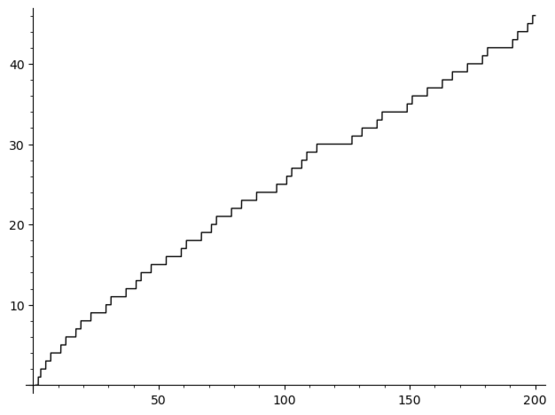}
\end{center}

\noindent
Riemann showed that the prime counting function can be written as a main term plus a sum of oscillating ``corrections", one for each nontrivial zero of the Riemann zeta function.
For example, when we include the first 13 correction terms we get a smooth function that 
approximates \(\pi(n)\) quite well for \(n \le 20\).

\vskip 1em
\begin{center}
\includegraphics[width=21em]{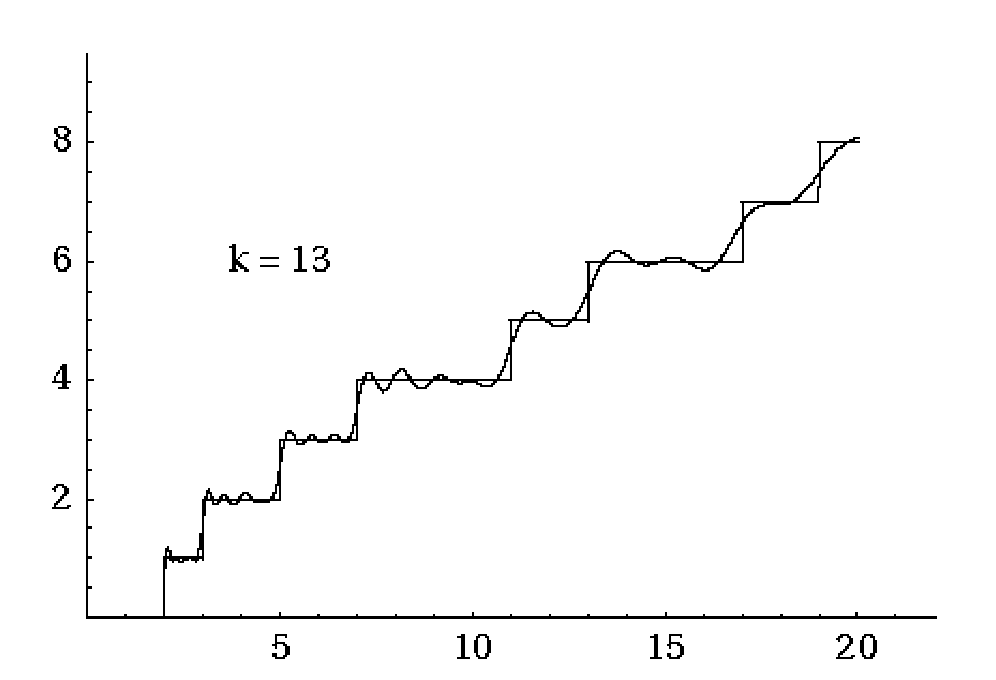}
\end{center}

\vskip 1em
When we include the first 118, we get a function that approximates \(\pi(n)\) fairly well
well even for \(n\) up to \(230\).

\vskip 1em
\begin{center}
\includegraphics[width=20em]{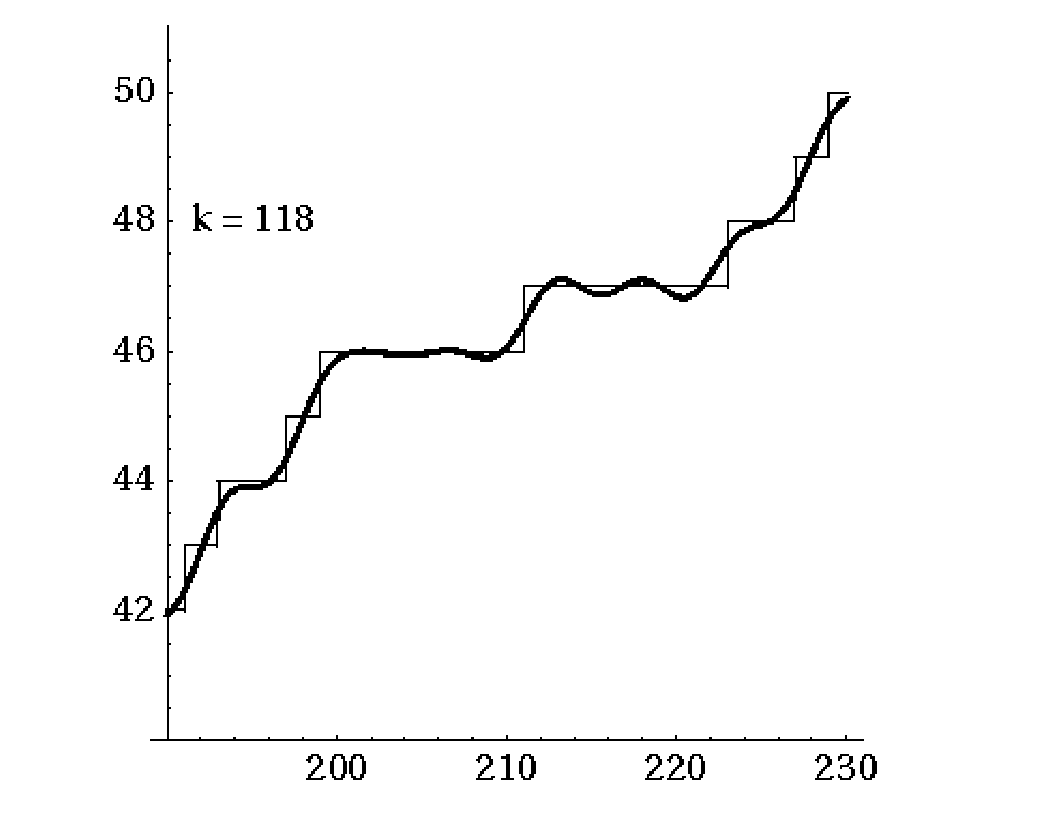}
\end{center}

\vskip 1em
If primes are like ``particles''---points where \(\pi(n)\) jumps discontinuously---then 
zeros of the Riemann zeta function correspond to ``waves".   This is not just poetry: there is a mathematical connection to the Fourier transformation, which underlies wave-particle duality in quantum mechanics.

Since the prime counting function \(\pi(n)\) equals this ``main term":
\[  \mathrm{li}(n) = \int_0^n \frac{dt}{\ln t} \]
plus corrections coming from the nontrivial Riemann zeta zeros,  
knowing the location of these zeros would give more information 
about the prime counting function.   Indeed, the Riemann Hypothesis:

\begin{center}
All nontrivial Riemann zeta zeros lie on the line \(\mathrm{Re}(z) = \frac{1}{2}\).
\end{center}
is equivalent to this claim:
\begin{center}
For some \(C > 0\)  and all \(n \ge 1\), \(| \pi(n) - \mathrm{li}(n) | \le C \sqrt{n} \ln n \).
\end{center}

More simply put, the Riemann Hypothesis says that the wavelike corrections to a 
simple approximation to the prime counting function are not very large.   This would have many implications for number theory: for example, it would imply a better bound on the biggest gaps between primes than is currently known.  There is thus is ample reason for wanting to prove this hypothesis. 

 So far, numerical computations far suggest that the Riemann Hypothesis is true: for example, the first 12 trillion nontrivial zeros of the Riemann zeta function lie on the line  \(\mathrm{Re}(z) = \frac{1}{2}\).   But this hypothesis has been stubbornly resistant to proof.  One reason is that we do not have a deep understanding of what the Riemann zeta zeros \emph{mean}---other than just providing oscillatory corrections to the prime counting function.

\section{The Weil Conjectures}

The Weil Conjectures are interesting because they include a variant of the Riemann Hypothesis that, while still difficult, has actually been proved.   In this variant, the count of solutions of some polynomial equations in several variables has a ``main term"  and some oscillatory ``correction terms".  One difference is that there are only finitely many correction terms.  Another is that we know more about the meaning of the terms: they come from things called ``motives".  
 
Let's look at an example: the polynomial equation 
\[  y^2 + y = x^3 + x .\]
To get a finite number of solutions, we let the variables \(x\) and \(y\) take values in a 
finite field.  All we need to know about finite fields is that there is one called \(\F_q\) of cardinality \(q = p^n\) for each prime \(p\) and each integer \(n = 1,2,3,\dots\).  Let's take \(p = 2\) and count the solutions of the above equation with \(x, y \in \F_q\).

\[ \begin{array}{rrr}
n & \;\;\;\; q = 2^n \! & \;\;\;\; \textrm{number of solutions} \\
1 & 2 & 4 \\
2 & 4 & 4 \\
3 & 8 & 4 \\
4 & 16 & 24 \\
5 & 32 & 24 \\
6 & 64 & 64 \\
7 & 128 & 144 \\
8 & 256 & 224 \\
9 & 512 & 544 \\
10 & 1024 & 1024 \\
11 & 2048 & 1984 \\
12 & 4096 & 4224 \\
\end{array}
\]

Since \(y^2 + y = x^3 + x \) is one equation with two unknowns, we might naively guess that in the field with \(p^n\) elements it has \(p^n\) solutions.   As you can see from the table, this is pretty close!   This approximation is the ``main term".  Let's subtract it off from the true number of solutions and get the ``correction term'':
\[ \begin{array}{rr}
n & \;\;\;\; \textrm{correction term} \\
1 & 2 \\
2 & 0 \\
3 & -4 \\
4 &  8 \\
5 & -8 \\
6 & 0 \\
7 & 16 \\
8 & -32 \\
9 & 32 \\
10 & 0 \\
\end{array}
\]
We can immediately notice some interesting things about this correction term.
First, its magnitude grows in a roughly exponential way while it oscillates between being positive and negative.  Second, when it's not zero it's a power of two.   Third, it starts out
nonzero.   All these patterns are real---and they suggest a correction term proportional to a cosine times an exponential, or in other words, proportional to 
\[      \alpha^n + \overline{\alpha}^n \]
for some complex number \(\alpha\).   With some experimentation we can guess \(\alpha = -1 + i\).  With this choice of \(\alpha\), here is what we get:
\[ \begin{array}{rrrr}
n & \qquad \qquad \alpha^n &\;\;\;\; \alpha^n + \overline{\alpha}^n \! & \;\;\;\; \textrm{correction term} \\
1 & -1 + i & -2 &  2 \\
2 & -2i & 0 & 0 \\
3 & 2 + 2i & 4 & -4 \\
4 & -4 & -8 & 8 \\
5 & 4 - 4i & 8 & -8 \\
6 & 8i & 0 & 0 \\
7 & -8 - 8i & -16 &  16 \\
8 & 16 & 32 & -32 \\
9 & -16 - 16i & -32 & 32 \\
10 & -32i & 0 & 0 \\
\end{array}
\]
So, it seems that the number of solutions of \(y^2 + y = x^3 + x \) in the field with \(2^n\) element is exactly
\[    2^n - \alpha^n - \overline{\alpha}^n \]
where \(\alpha = -1+i\).   Indeed this is true!  

What is special about the equation \(y^2 + y = x^3 + x\)?   Crucially, it gives an example of an ``elliptic curve''.  The precise definition of an elliptic curve is a bit technical, but the rough idea is that a polynomial equation with integer coefficients in two variables gives an elliptic curve if its space of \emph{complex} solutions is a torus with one point removed.   We then say the elliptic curve ``over \(\C\)'' is the torus.  In other words, it's the set of complex solutions with an extra point included.

Elliptic curves are not hard to find.   For example, most cubic equations in two variables give elliptic curves, such as
\[   y^2 = x^3 + ax + b \]
for any \(a,b \in \Z\).   In 1933, Helmut Hasse proved a theorem that vastly generalizes the result we have empirically observed:

\begin{thm}[\textbf{Hasse's Theorem}]  Given a polynomial equation with integer coefficients in two variables that gives an elliptic curve, for any prime power \(q = p^n\),
the number of solutions in \(\F_q\) of this equation is
\[  p^n - \alpha^n - \overline{\alpha}^{\, n} \]
where \(\alpha \in \C\) has \(|\alpha| = \sqrt{p}\).  
\end{thm}

In this situation we say the elliptic curve ``over \(\F_q\)'' is the set of solutions of the polynomial equation in \(\F_q\) together with one extra point.  So, the number of points of the elliptic curve over \(\F_q\) is
\[   p^n - \alpha^n - \overline{\alpha}^{\, n} + 1. \]

Now comes something amazing.  The four terms in the above formula correspond, in some subtle and mysterious way, to four pieces of the elliptic curve over \(\C\), which is a torus:
\begin{center}
\includegraphics[width=10em]{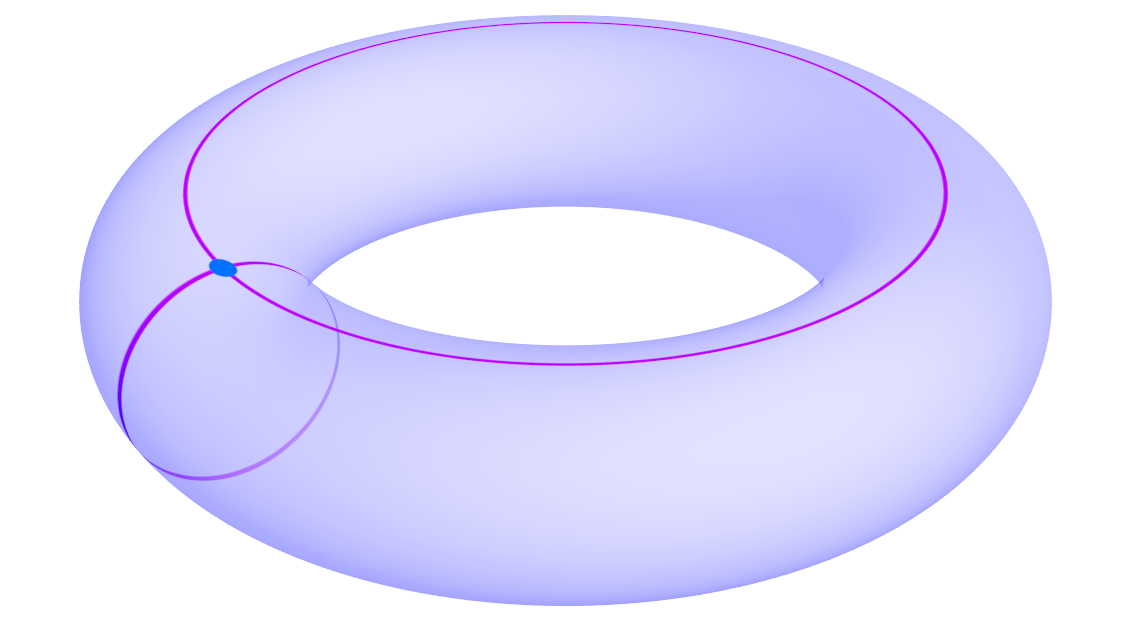}
\end{center}
\noindent 
These four pieces are the point, the two circles with that point removed, and the whole torus with the two circles and point removed.  The pieces of dimension \(k\) correspond to terms in the above formula that grow like \(p^{nk/2} = q^{k/2}\).   

In a rough, handwaving way it may seem plausible that the piece of dimension \(2\) contributes \(q\) points to the elliptic curve.  After all, in the complex version of the picture, the torus minus the point and two circles can be identified with \(\C\).  Perhaps by analogy when we work over the field \(\F_q\) the elliptic curve will have a similar piece that is a copy of \(\F_q\), containing \(q\) points.    It may also seem believable that the elliptic curve over \(\F_q\) should have a piece of dimension \(0\) containing one point: namely, the extra point we deliberately added on.  

The mystery is how the curve over \(\F_q\) could have two pieces corresponding to the oscillating correction terms \(-\alpha^n - \overline{\alpha}^{\, n}\) in the count of points.  After all, these two terms are \emph{complex}, and while they sum to an integer, this integer can be \emph{negative}.   So, we cannot think of these terms as literally counting points in some subsets of the elliptic curve over \(\F_q\).

In short, the individual terms in the formula that counts points of an elliptic curve over 
\(\F_q\) do not simply count points in subsets of this curve.   Grothendieck's idea was that they correspond to subtler pieces of the curve, called \emph{motives}.    Bizarrely, these pieces can have a negative number of points.

What are these things called motives?   What could they possibly be?  Since they are supposed to explain oscillating terms in the count of points, we might hope they are connected to vector spaces---somewhat like how quantum mechanics uses vector spaces to describe oscillating states of physical systems.  This turns out to be true!  

But motives are also connected to topology---algebraic topology, to be precise.  And this is no coincidence: algebraic topology lets us turn topological spaces into vector spaces. To any topological space \(X\) we  can associate a sequence of vector spaces \(H^k(X)\) called its ``cohomology groups'', which keep track of \(k\)-dimensional aspects of its topology.  

Consider in particular an elliptic curve over \(\C\): that is, a torus.  Its 0th cohomology group is a 1-dimensional vector space, and we can get a basis vector for this space from any point in the torus.  Its 1st cohomology group is a 2-dimensional vector space, and we can get basis vectors for this from the two circles shown above---though many other choices of circles would work equally well.  Finally, its 2nd cohomology group is a 1-dimensional vector space, and a basis for this is given by a 2-dimensional surface, namely the torus itself.

Of course, the mysterious part is that the topology of the elliptic curve over \(\C\) is relevant to its count of points over a finite field!  This is very hard to understand.  But this turned out to be the key to generalizing Hasse's theorem.  

First, with a lot of work from roughly 1940 to 1948, Andr\'e Weil generalized Hasse's result to algebraic curves of arbitrary genus \(g\): that is, polynomial equations in two variables whose space of complex solutions, plus one extra point, look like a smooth surface with \(g\) handles.   Here is an example with \(g=2\):

\begin{center}
\includegraphics[width=12em]{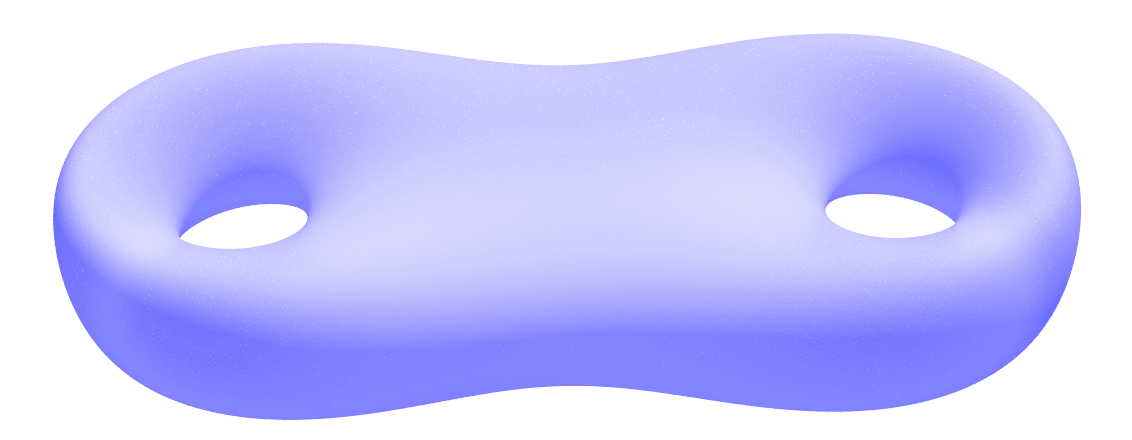}
\end{center}

The 0th and 2nd cohomology groups of such a curve are still 1-dimensional, but the the 1st cohomology group has dimension \(2g\).  So, when we count points of an algebraic curve defined over a finite field, we expect more terms in the formula.   Weil proved that this is true:

\begin{thm}[\textbf{Weil's Theorem}]
Given a polynomial equation with integer coefficients in two variables that gives an algebraic curve of genus \(g\), for any prime power \(q = p^n\) the number of points of this curve 
over \(\F_q\) is
\[   p^n - \alpha_1^n - \cdots - \alpha_{2g}^n + 1 \]
where all the numbers \(\alpha_i \in \C\) have \(|\alpha_i| = \sqrt{p}\).   
\end{thm}

The simplest example is the curve of genus zero.  This is called the ``projective line'' \(\P^1\), and over the complex numbers it is usually called the ``Riemann sphere'', since it is a sphere consisting of a copy of \(\C\) together with a point at infinity:

\begin{center}
\includegraphics[width=9em]{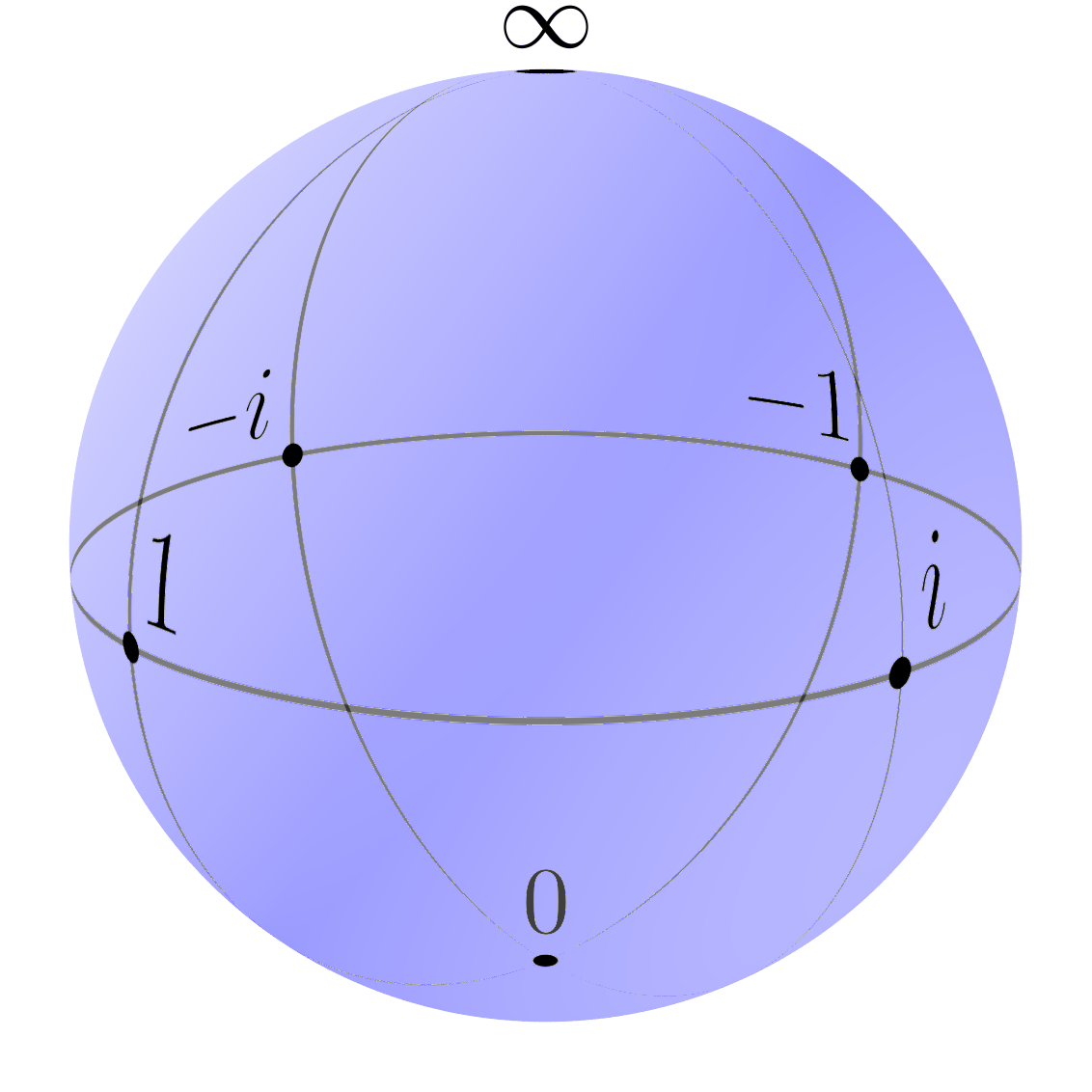}
\end{center}

\noindent 
Similarly over \(\F_q\) we can chop the projective line into two subsets---a copy of \(\F_q\) and a point at infinity---so it has
\[     q + 1 \]
points.  Thus, in this particular case the two ``motives'' making up the projective
line can be loosely thought of as as actual subsets.  But this is deceptive: there is nothing special about the point at infinity here; any other point would do just as well.

This trivial example hints at a higher-dimensional generalization of Weil's
theorem.  Besides the projective line, there is also a ``projective plane'',
which consists of a plane together with a projective line at infinity.   
Over the real numbers the plane is \(\R^2\), but over
the complex numbers it is \(\C^2\), while the projective
line is \(\C + 1\): that is, the disjoint union of \(\C\) and a point.  So,
we can write the complex projective plane \(\C\P^2\) as
\[    \C\P^2 = \C^2 + \C + 1 .\]
Similarly we can chop the projective plane over the finite field \(\F_q\), into three
pieces: the plane \(\F_q^2\), which 
has \(q^2\) points, and the projective line \(\F_q + 1\):
\[   \F_q\P^2 = \F_q^2 + \F_q + 1 .\]
Thus, \(\F_q\P^2\) has 
\[  q^2 + q + 1 \]
points.   

Higher-dimensional projective spaces follow the same pattern.   But elliptic curves illustrate the subtleties we must confront in generalizing Weil's result to higher dimensions.  As we've seen, elliptic curves have both even- and odd-dimensional cohomology over \(\C\).
For this reason, over \(\F_q\) their number of points is not merely polynomial in \(q\), but
also has terms with magnitude proportional to powers of \(\sqrt{q}\).

To generalize his result for curves to higher dimensions, Weil needed to consider algebraic varieties.  Given a collection of polynomial equations with integer coefficients in \(n\) variables, we can study their solutions in---or as the experts say, ``over''---any field \(\F\).   The solutions form a subset of \(\F^n\), which is called an ``affine algebraic variety'' over \(\F\).   However, it is often useful to include additional points at infinity, as we have already seen in some examples.  If we do this correctly, we get a ``projective algebraic variety''.  

If we work over the field \(\F = \C\), and our projective algebraic variety \(X\) is smooth in a suitable sense, then it is a \(2d\)-dimensional manifold, and it can have nontrivial cohomology groups \(H^k(X)\) only for \(0 \le k \le 2d\).    But an algebraic variety over a finite field \(\F_q\) has a finite set of points, which we can count.   

This line of thought led Weil to conjecture a grand generalization of his previous result for curves.   He made this conjecture in 1949.   Building on decades of work by Grothendieck and others, the proof was finally completed by Deligne in 1974.

\begin{thm}[\textbf{Riemann Hypothesis for Varieties over Finite Fields}]
Given a collection of polynomial equations with integer coefficients defining a smooth projective variety, for any prime power \(q = p^n\), the number of points of this variety over \(\F_q\) is
\[ \sum_{k = 0}^{2d} \sum_{i = 1}^{\beta_k} \; (-1)^k \alpha_{ik}^n \]
where \(|\alpha_{ik}| = p^{k/2}\) and \(\beta_k\) is the dimension of the \(k\)th cohomology group of the corresponding projective variety over \(\C\).
\end{thm}

This theorem is often phrased in terms of a ``zeta function'' associated to the variety \(X\), a relative of the function Riemann was interested in.  This zeta function is essentially a trick for keeping track of the count of points on \(X\).    Theorem 3 can be restated as saying that this zeta function has zeros and poles on the lines \(\textrm{Re}(z) = \frac{k}{2}\): zeros when \(k\) is odd and poles when \(k\) is even.  

However, the connection to motives is most easily visible if we inquire about the terms in the above sum.    Grothendieck's dream was that we can somehow break \(X\) into pieces of dimension \(k = 0, 1, \dots, 2d\), with the \(k\)-dimensional pieces contributing all the terms of the form \((-1)^k \, \alpha_{ik}^n\) to the number of points.   As we've seen, these pieces cannot in general be subsets of \(X\).  They must be something else: motives.

\section{Motives}

But what is a motive?   Here we need to turn up the heat and assume some familiarity 
with category theory.   For any prime power \(q\) there is a category \(\Var\) of smooth projective varieties over \(\F_q\) and the usual maps between these.  Starting from \(\Var\) we can construct a category  \(\Mot\), called the category of 
``pure motives", and a functor
\[ h \colon \Var \to \Mot^{\rm{op}},  \]
where the use of the opposite category is purely traditional.   This category \(\Mot\), and also its opposite, have many features resembling the category of vector spaces and linear maps:
\begin{itemize}
\item
They  are ``linear categories": the hom-sets are vector spaces, and composition is bilinear. 
\item
They are ``Karoubian'' or ``Cauchy complete": they have direct sums, and any \(\pi \colon X \to X\) with \(\pi^2 = \pi\) is projection onto \(Y\) for some direct sum decomposition \(X \cong Y \oplus Z\).
\item
They are ``symmetric monoidal": they have a well-behaved tensor product \(\otimes\), coming from the cartesian product of smooth projective varieties.
\end{itemize}
\noindent
Thus, we should think of the functor \(h\) as taking us from the world of smooth projective varieties to the world of linear algebra.  

To get a sense for how Cauchy completeness works, note that 
we can map all of \(\P^1\) to a single point \(p \in \P^1\), defining 
\[   \pi \colon \P^1 \to \P^1  \]
by \(\pi(x) = p\) for all \(x \in \P^1\).   This map clearly has \(\pi^2 = \pi\), so 
the morphism \(h(\pi) \maps h(\P^1) \to h(\P^1)\) has \(h(\pi)^2 = h(\pi)\).
By Cauchy completeness, \(h(\pi)\) is the projection onto 
some summand in a direct sum decomposition of \(h(\P^1\)).  With a bit of category
theory we can check that this summand must be isomorphic to \(h(1)\), where \(1 \in \Var\) is the 1-point variety.  So,
\[  h(\P^1) = h(1) \oplus \L  \]
where \(\L\) is some other motive.

This other motive \(\L\) is called the ``Lefschetz motive''.  It is not \(h\) of any 
smooth projective variety!   We can loosely visualize it as the projective line with a point removed---and thus, an ordinary line.  Earlier, we chopped the projective line
into a point and an ordinary line.  The latter was not a projective variety---but the decomposition \(h(\P^1) = h(1) \oplus \L\) is valid in the category
of motives. 

With more work, we can show
\[   h(\P^n) \; \cong \; h(1) \; \oplus \; \L \; \oplus \; \L^{\otimes 2} \; \oplus \; \cdots \; \oplus \; \L^{\otimes n}  \]
and this corresponds to the formula for the number of points of \(\P^n\) over \(\F_q\):
\[  1 \; + \; q \; + \; \cdots \; + \; q^n .\]
But curves and other varieties typically decompose into motives 
that are \emph{not} just tensor powers of the Lefschetz motive \(\L\), because
their number of points over \(\F_q\) is not just a polynomial in \(q\).

Grothendieck showed the Riemann Hypothesis for finite fields would follow from the so-called ``Standard Conjectures".   Among other things, these conjectures would imply:

\begin{itemize}
\item
Every variety \(X\) has \(h(X) \cong X_0 \oplus \cdots \oplus X_n\) where the motive \(X_k\) has dimension \(k\), or technically speaking ``weight" \(k\), meaning that it contributes terms proportional to \(\alpha^n\) with \(|\alpha| = p^{k/2}\) to the count of points of \(X\) over \(\F_{p^n}\).    
\vskip 1em
\item
The category \(\Mot\) is ``abelian": it has well-behaved kernels, cokernels, subobjects and quotient objects.
\vskip 1em 
\item
The category \(\Mot\) is ``semisimple": 
every motive is a finite direct sum of so-called ``simple'' motives that have only two subobjects, 0 and that motive itself.  Thus, each motive \(X_k\) above can be further decomposed into a direct sum of simple motives.
\end{itemize}

If the Standard Conjectures are true, we can take any \(d\)-dimensional smooth projective variety \(X\) over the algebraic completion of \(\F_p\) and compute its number of points over each field \(\F_{p^n}\) as follows.   First, we break the motive \(h(X)\) into a direct sum of motives \(X_k\) as above.  Each of these comes equipped with a special morphism
\[     F_k \maps X_k \to X_k \]
called the ``Frobenius''.   The number of points of \(X\) over \(\F_{p^n}\) can then be
expressed in terms of the \(n\)th power of the Frobenius as follows:
\[    \sum_{k = 0}^{2d} (-1)^n \mathrm{tr}(F_k^n)\]
where the trace is defined by carrying ideas from linear algebra to the category of pure motives.   This explains the exponentially growing yet also perhaps oscillating 
terms in the Riemann Hypothesis for varieties over finite fields (see Theorem 3).

Another consequence of the Standard Conjectures is that \(h\) is the ``universal Weil cohomology theory'': that is, the maximally informative cohomology theory for smooth projective varieties obeying a certain list of axioms.   

Alas, Grothendieck was unable to prove the Standard Conjectures.   It is still not known if 
\(\Mot\) has the desirable properties just listed.   In 1974 Deligne proved the Riemann Hypothesis for varieties over finite fields in a way that sidestepped these questions.

Thus, motives remain deeply mysterious---yet the \emph{idea} of motives, or what Grothendieck called the ``yoga'' of motives, has been very powerful in the hands of
skilled mathematicians, from Deligne to Voevodsky and others.   There are also bold visions that go even further.   For example, in 1995 Yuri Manin outlined a dream of proving the actual Riemann Hypothesis by generalizing motives over fields to motives over a poorly-understood entity called the ``field with one element", which is not really a field.  More recently Alain Connes has been trying to prove the Riemann Hypothesis using these ideas, combined with ideas taken from quantum physics. 

So, motives continue to tantalize and inspire!

\section{Further reading}

Here we give references to works mentioned above and also some suggestions for further
reading---mainly on pure motives, since mixed motives, motivic Galois groups and
motivic cohomology go beyond the limited scope of this article.

\subsection{Introduction}

The Grothendieck quote about motives is from Lisker's translation of Grothendieck's 1500-page memoir \textsl{R\'ecoltes et Semailles}, which has many interesting discussions of motives:

\begin{itemize}
\item Alexander Grothendieck, \textsl{Récoltes et Semailles}, 1986.  Available at
\href{https://agrothendieck.github.io/divers/ReS.pdf}{https://} \href{https://agrothendieck.github.io/divers/ReS.pdf}{agrothendieck.github.io/divers/ReS.pdf}. Partial translation into English by Roy Lisker available at \href{http://matematicas.unex.es/~navarro/res/lisker1.pdf}{http://matematicas.unex.es/$\sim$navarro/res/lisker1.pdf}. %\href{http://matematicas.unex.es/~navarro/res/lisker1.pdf}{lisker1.pdf}.
\end{itemize}

\subsection{The Riemann Hypothesis}

Riemann propounded his famous hypothesis here:

\begin{itemize}
\item
Bernard Riemann, Ueber die Anzahl der Primzahlen unter einer gegebenen Gr\"osse, \textsl{Monatsberichte der Berliner Akademie}, November 1859.  Available in German and English at \href{https://www.maths.tcd.ie/pub/HistMath/People/Riemann/Zeta/}{https://www.maths.tcd.ie/pub/HistMath/People/}  \href{https://www.maths.tcd.ie/pub/HistMath/People/Riemann/Zeta/}{Riemann/Zeta/}.
\end{itemize}
This is an enjoyable introduction to the Riemann Hypothesis and its connection to the prime counting function:
\begin{itemize}
\item Barry Mazur and William Stein,
\textsl{Prime Numbers and the Riemann Hypothesis}, Cambridge U.\ Press,
Cambridge, 2016.
\end{itemize}
The picture of the prime counting function was drawn for this pape by Dima Pasechnik.
The two pictures of approximations to the prime counting function are still shots of animations that are available online:
\begin{itemize}
\item Dan Rockmore, Chance in the primes, 
\href{https://chance.dartmouth.edu/chance\_news/recent\_news/chance\_news\_10.10.html}{https://chance.dartmouth.edu/chance} \href{https://chance.dartmouth.edu/chance\_news/recent\_news/chance\_news\_10.10.html}{\_news/recent\_news/chance\_news\_10.10.html}.
\end{itemize}

\subsection{The Weil Conjectures}

James Milne has written some very useful papers on motives and the Weil Conjectures.
This is an insightful history of the latter:

\begin{itemize}
\item James Milne, The Riemann Hypothesis over finite fields: from Weil to the present day,
in \textsl{The Legacy of Bernhard Riemann after One Hundred and Fifty Years}, vol.\ II, eds.\
Lizhen Ji, Frans Oort and Shing-Tung Yau, International Press, Somerville, Massachusetts,  
 2015, pp.\ 487--565.  Available at \href{https://www.jmilne.org/math/xnotes/pRH.html}{https://www.jmilne.org/math/xnotes/}\href{https://www.jmilne.org/math/xnotes/pRH.html}{pRH.html}. 
\end{itemize}
Instead of providing references to the work of Hasse, Weil, Grothendieck, Deligne and others connected to the Weil Conjectures, I defer to this work.

The pictures of an elliptic curve, genus-2 curve and Riemann sphere were created for this paper by Simon Burton.

\subsection{Motives}

This is an excellent quick introduction to the definition of pure motives and the Standard
Conjectures:

\begin{itemize}
\item James Milne, Motives: Grothendieck's dream, in \textsl{Open Problems and 
Surveys of Contemporary Mathematics}, eds.\ Lizhen Ji, Yat-Sun Poon and Shing-Tung Yau, International Press, Somerville, Massachusetts, 2013, pp.\
325--342.  Available at \href{https://www.jmilne.org/math/xnotes/mot.html}{https://www.jmilne.org/math/xnotes/mot.html}.
\end{itemize}
Milne has a real flair for saying only what absolutely needs to be said
to get the main ideas across.  However, like most treatments of motives, his introduction assumes a good understanding of algebraic geometry and cohomology theory.  For more try this excellent book, which has roughly similar prerequisites, but fills in more details and goes further:
\begin{itemize}
\item Jacob P.\ Murre, Jan Nagel and Chris A.\ M.\ Peters, \textsl{Lectures on the Theory of Pure Motives}, University Lecture Series 61, AMS, Providence, 2013.
\end{itemize}
To dig deeper, try this:

\begin{itemize}
\item Yves Andr\'e, \textsl{Une Introduction aux Motifs (Motifs Purs, Motifs Mixtes, 
P\'eriodes)}, \textsl{Panoramas et Synth\`eses} \textbf{17}, Soci\`et\`e Math.\ de France, Paris, 2004.
\end{itemize}

There are various kinds of pure motives arising from various equivalence relations on algebraic cycles.  In the body of this paper I was implicitly defining them using ``homological equivalence''.   If instead one defines them using a coarser equivalence relation called ``numerical equivalence'', one can prove the resulting category of pure motives is abelian
and semi-simple.  This was done by Jannsen:
\begin{itemize}
\item Uwe Jannsen, Motives, numerical equivalence, and semi-simplicity, \textsl{Invent.\ Math.} \textbf{107} (1992), 447--452.  Available at \href{https://eudml.org/doc/143974}{https://eudml.org/doc/143974}.
\end{itemize}
Furthermore, one of the unproven Standard Conjectures, called Conjecture D, asserts that homological and numerical equivalence actually agree!  Thus, it is very interesting to study pure motives defined using numerical equivalence, which can be understood quite explicitly, at least over finite fields:

\begin{itemize}
\item 
James Milne, Motives over finite fields, \textsl{Proc.\ Sympos.\ Pure Math.} \textbf{55},
Part 1, AMS, Providence, 1994, pp.\ 401--459.  Available at 
\href{https://www.jmilne.org/math/articles/1994a.html}{https://www.jmilne.} \href{https://www.jmilne.org/math/articles/1994a.html}{org/math/articles/1994a.html}.
\end{itemize}

Manin's dream of applying motivic ideas to the actual Riemann Hypothesis is explained
here:
\begin{itemize}
\item Yuri Manin, Lectures on zeta functions and motives, \textsl{Ast\'erisque} \textbf{228} (1995), 121--163.  Available at 
\href{http://www.numdam.org/item/?id=AST_1995__228__121_0}{http://www.numdam.org/item/?id=AST\_1995\_\_228}
\href{http://www.numdam.org/item/?id=AST_1995__228__121_0}{http://www.numdam.org/item/?id=AST\_}{\_\_121\_0}
\end{itemize}
Connes has been working with a number of collaborators to prove the Riemann
Hypothesis using ideas connected to motives.  An early report is his book with Marcolli:
\begin{itemize}
\item Alain Connes and Matilde Marcolli, \textsl{Noncommutative Geometry, Quantum Field Theory and Motives}, AMS, Providence, 2007.  Available at \href{https://alainconnes.org/wp-content/uploads/bookwebfinal-2.pdf}{https://alainconnes.org/wp-content/uploads/bookwebfinal-2.pdf}.
\end{itemize}

\vskip 2em
\subsection*{Acknowledgements}

I thanks Peter Jipsen, Alexander Kurz, Andrew Mosher, Marco Panza, Ahmed Sebbar 
and Daniele Struppa for inviting me to speak on this topic at the Grothendieck
Conference at Chapman University on May 25th, 2022.  I thank Simon Burton and Dima Pasechnik for figures, and David Egolf and the anonymous referee for pointing out some
explanations that needed clarification.

\vfill

\end{document}